\documentclass[10pt]{article} 
\usepackage{amsfonts} 
 
\usepackage{amscd}

\newtheorem{theorem}{Theorem}[section]  
  
\newtheorem{lemma}[theorem]{Lemma}  
\newtheorem{proposition}[theorem]{Proposition}  
  
\newtheorem{definition}{Definition}[section]  
\newtheorem{remark}{Remark}[section]  
\begin{document}
\title{Affine structures on abelian Lie Groups}
\author{Elisabeth Remm - Michel Goze \\
Universit\'{e} de Haute Alsace, F.S.T.\\ 
4, rue des Fr\`{e}res Lumi\`{e}re - 68093 MULHOUSE - France}
\date{}
\maketitle

The Nagano-Yagi-Goldmann theorem states that on the torus $\mathbb{T}^{2},$
every affine (or projective) structure is invariant or is constructed on the
basis of some Goldmann rings [12]. It shows the interest to study the
invariant affine structure on the torus $\mathbb{T}^{2}$ or on abelian Lie
groups. Recently, the works of Kim [8] and Dekimpe-Ongenae [5] precise the
number of non equivalent invariant affine structures on an abelian Lie group
in the case these structures are complete. In this paper we propose a
study of complete and non complete affine structure on abelian Lie groups
based on the geometry of the algebraic variety of finite dimensional
associative algebras.

\section{Affine structure on Lie groups and Lie algebras}

The Lie group $Aff(\mathbb{R}^{n})$ is the group of affine transformations
of $\mathbb{R}^{n}.$ It is constituted of matrices 
\[
\left( 
\begin{array}{cc}
A & b \\ 
0 & 1
\end{array}
\right)
\]
with $A\in GL(n,\mathbb{R}),$ $b\in \mathbb{R}^{n}.$ It acts on the real
affine space $\widetilde{\mathbb{R}}^{n}$ by 
\[
\left( 
\begin{array}{cc}
A & b \\ 
0 & 1
\end{array}
\right) \left( 
\begin{array}{c}
v \\ 
1
\end{array}
\right) =\left( 
\begin{array}{c}
Av+b \\ 
1
\end{array}
\right)
\]
where $(v,1)^{t}\in \widetilde{\mathbb{R}}^{n+1}.$

Its Lie algebra, noted $aff(\mathbb{R}^{n}),$ is the linear algebra 
\[
aff(\mathbb{R}^{n})=\left\{ \left( 
\begin{array}{cc}
A & b \\ 
0 & 0
\end{array}
\right) ;A\in gl(n,\mathbb{R)},b\in \mathbb{R}^{n}\right\}
\]

\begin{definition}
An affine structure on a Lie algebra $\frak{g}$ \ is a morphism 
$$
\Psi :\frak{g}\mathcal{\rightarrow }aff(\mathbb{R}^{n})
$$
of Lie algebras.
\end{definition}

\noindent \textbf{Remark. }Let us consider an affine representation of a Lie group $G$
 that is an homomorphism $\varphi :G\mathcal{\rightarrow }Aff(\mathbb{R}%
^{n}).$ For every $g\in G$ , $\varphi (g)$ is an affine transformation on
the affine space $\widetilde{\mathbb{R}}^{n+1}.$ This representation induces
an affine structure on the Lie algebra $\frak{g}$ of $G$.

\begin{proposition}
The Lie algebra $\frak{g}$ is provided with an affine structure if and only
if the underlying vector space $A(\frak{g}\mathcal{)}$ is a
left symmetric algebra, that is there exists a bilinear mapping 
\[
\begin{array}{ccc}
A(\frak{g}\mathcal{)\times }A(\frak{g}\mathcal{)} & \rightarrow & A(\frak{g}%
\mathcal{)} \\ 
\left( X,Y\right) & \mapsto & X.Y
\end{array}
\]
satisfying

1) $X.\left( Y.Z\right) -Y.\left( X.Z\right) =\left( X.Y\right) .Z-\left(
Y.X\right) .Z$

2) $X.Y-Y.X=\left[ X,Y\right] $

for every $X,Y,Z\in A(\frak{g}\mathcal{)}$
\end{proposition}

If $\Psi $ is a morphism giving an affine structure on $\frak{g}$, then the
left symmetric product on $A(\frak{g}\mathcal{)}$ is defined as this: 
$$
\forall X\in \frak{g},\Psi \left( X\right) =\left( 
\begin{array}{cc}
A\left( X\right) & b\left( X\right) \\ 
0 & 0
\end{array}
\right)
$$
and we put 
$$
X.Y=b^{-1}\left( A\left( X\right) .b\left( X\right) \right)
$$
where $b:A(\frak{g}\mathcal{)\rightarrow }\mathbb{R}^{n}$ is supposed
to be bijective.

The fact that $\Psi $ is a representation implies that $X.Y$ is a left
symmetric product. Conversely, if $X.Y$ is a left symmetric product on $A(%
\frak{g}\mathcal{)},$ and if $L_{X}$ indicates the left representation : $%
L_{X}\left( Y\right) =X.Y,$ then the map 
\[
X\rightarrow \left( 
\begin{array}{cc}
L_{X} & X \\ 
0 & 0
\end{array}
\right)
\]
defines an affine structure on the Lie algebra $\frak{g}.$

If $\Psi $ is an affine structure on $\frak{g}$, it defines a representation 
$$
\Psi \left( X\right) =\left( 
\begin{array}{cc}
A\left( X\right) & b\left( X\right) \\ 
0 & 0
\end{array}
\right)
$$
and a left symmetric product $X.Y$. This last induces an affine
representation 
\[
\left( 
\begin{array}{cc}
L_{X} & X \\ 
0 & 0
\end{array}
\right)
\]
which is equivalent to $\Psi $ (and equal if $b=Id$).

\begin{definition}
An affine structure on $\frak{g}$ is called complete if the endomorphism 
\[ 
\begin{array}{cccc}
\theta _{X}: & A(\frak{g}\mathcal{)} & \rightarrow & A(\frak{g}\mathcal{)}
\\ 
& Y & \mapsto & Y+Y.X
\end{array}
\]
is bijective for every $X\in A(\frak{g}\mathcal{)}.$
\end{definition}

This is equivalent to one of the following properties .

a) 
\[
\begin{array}{cccc}
R_{X}: & A(\frak{g}\mathcal{)} & \rightarrow & A(\frak{g}\mathcal{)} \\ 
& Y & \mapsto & Y.X
\end{array}
\]
is nilpotent for all $X\in A(\frak{g}\mathcal{)}$

b) $tr(R_{X})=0$ for all $X$ $\in A(\frak{g}\mathcal{)}.$

\medskip

\noindent \textbf{Remarks }

1. In [13] we have given examples of non complete affine structures on some
nilpotent Lie algebras of maximal class. In particular, on the 3-dimensional
Heisenberg algebra  the affine structure associated to the following
representation 
\[
\left( 
\begin{array}{llll}
a(x_{1}+x_{2}) & a(x_{1}+x_{2}) & 0 & x_{1} \\ 
a(x_{1}+x_{2}) & a(x_{1}+x_{2}) & 0 & x_{2} \\ 
\alpha x_{1}+(\beta -1)x_{2} & \beta x_{1}+(\alpha +1)x_{2} & 0 & x_{3} \\ 
0 & 0 & 0 & 0
\end{array}
\right)
\]
is non complete.

2. The complete affine structure on $\frak{g}$ corresponds to the
simply-transitive affine action of the connected corresponding Lie group $G.$

\section{Affine structure on abelian Lie algebras}

1. Let $\frak{g}$ be a real abelian Lie algebra. If $\frak{g}$ is provided
with an affine structure then the left symmetric algebra $A(\frak{g}\mathcal{%
)}$ is a commutative and associative real algebra. In fact we have 
\[
X.Y-Y.X=0=\left[ X,Y\right]
\]
and 
\[
X.\left( Y.Z\right) -Y.\left( X.Z\right) =\left( \left[ X,Y\right] .Z\right)
=0.
\]
This gives 
\[
X.\left( Z.Y\right) =\left( X.Z\right) .Y
\]
and $A(\frak{g}\mathcal{)}$ is associative.

Let $\Psi _{1}$ and $\Psi _{2}$ be two affine structures on $\frak{g}$. They
are affinely equivalent if and only if the corresponding commutative and
associative algebras are isomorphic. Thus the classification of affine
structures on abelian Lie algebras corresponds to the classification of
commutative and associative (unitary or not) algebras. If the affine
structure is complete, the endomorphisms $R_{X}$ are nilpotent. As $A(\frak{g%
}\mathcal{)}$ is also commutative, it is a nilpotent
associative commutative algebra. The classification of complete affine
structures on abelian Lie algebras corresponds to the classification of
nilpotent associative algebras. In this frame, we can cite the works of
Gabriel [7] and Mazzola [10] who study the varieties of unitary associative
complex laws and give their classification for dimensions less than 5, as
well as the works of Makhlouf and Cibils who study the deformations of these
laws and
 propose classifications for nilpotent associative complex algebras
[11].
 
\section{ Rigid affine structures}

\subsection{Definition}

Let us consider a fixed basis $\left\{ e_{1},...,e_{n}\right\} $ of the
vector space $\mathbb{R}^{n}.$ An associative law on $\mathbb{R}^{n}$ is
given by a bilinear mapping 
\[
\mu :\mathbb{R}^{n}\times \mathbb{R}^{n}\rightarrow \mathbb{R}^{n}
\]
satisfying $\mu (e_{i},\mu (e_{j},e_{k}))=\mu (\mu (e_{i},e_{j}),e_{k}).$ If
we put $\mu (e_{i},e_{j})=\sum C_{ij}^{k}e_{k}$, then the structural
constants $C_{ij}^{k}$ satisfy 
\[
(1)\quad {\sum_{l} }C_{il}^{s}C_{jk}^{l}-C_{ij}^{l}C_{lk}^{s}=0,%
\qquad s=1,...,n.
\]
Moreover if $\mu $ is commutative, we have 
\[
(2)\quad C_{ij}^{k}=C_{ji}^{k}
\]
Thus the set of associative laws on $\mathbb{R}^{n}$ is identified to the real
algebraic set embedded in $\mathbb{R}^{n^{3}},$ defined by the
polynomial equations (1) and (2). We note this set $\mathcal{A}^{c}(n).$

The law $\mu $ is unitary if there exists an $e\in \mathbb{R}^{n}$ such that $\mu
(e,x)=x.$ The set of unitary laws of $\mathcal{A}^{c}(n)$ is  noted by
$\mathcal{A}_{1}^{c}(n).$

The linear group $GL(n,\mathbb{R)}$ acts on $\mathcal{A}^{c}(n):$%
$$
\begin{array}{ccc}
GL(n,\mathbb{R)\times } \mathcal{A}^{c}(n) & \rightarrow & \mathcal{A}^{c}(n) \\ 
\left( f,\mu \right) & \mapsto & \mu _{f}
\end{array}
$$
where $\mu _{f}(e_{i},e_{j})=f^{-1}\mu (f(e_{i}),f(e_{j})).$

We note by $\theta (\mu )$ the orbit of $\mu $ by this action. The orbit is
isomorphic to the homogeneous space $\frac{GL(n,\mathbb{R)}}{G_{\mu }}$,
where $G_{\mu }=\left\{ f\in GL(n,\mathbb{R)}\quad /\quad \mu _{f}=\mu
\right\} .$ The topology of $\mathcal{A}^{c}(n)$ is the induced topology of $%
\mathbb{R}^{n^{3}}.$

\begin{definition}
The law $\mu \in \mathcal{A}^{c}(n)$ is rigid if $\theta (\mu )$ is open in $%
\mathcal{A}^{c}(n).$
\end{definition}

Let $\mu $ be a real associative algebra and let us note by $\mu _{\mathbb{C}%
}$ the corresponding complex associative algebra. If $\mu $ is rigid in $%
\mathcal{A}^{c}(n)$ then either $\mu _{\mathbb{C}}$ is rigid in the scheme $%
Ass_{n}$ of complex associative law, or $\mu _{\mathbb{C}}$ admits a
deformation $\widetilde{\mu }_{\mathbb{C}}$ which is never the
complexification of a real associative algebra.

This topological approach of the variety of associative algebras allows to
introduce the notion of rigidity on the affine structures.

\begin{definition}
An affine structure $\Psi $ on an abelian Lie algebra $\frak{g}$ is called
rigid if the corresponding associative algebra $A(\frak{g})$ is rigid in $%
\mathcal{A}^{c}(n).$
\end{definition}

\subsection{Cohomological approach}

It is well known that a sufficient condition for an associative algebra $%
\frak{a}$ to be rigid in $\mathcal{A}(n)$ is $H^{2}(\frak{a,a)=}0$, where $%
H^{\ast }(\frak{a,a)}$ is the Hochschild cohomology of $\frak{a.}$ Suppose
that $\frak{a}$ is commutative.\ If $\mu $ is the bilinear mapping defining
the product of $\frak{a,}$ then $\frak{a}$ is rigid if every deformation $%
\mu _{t}=$ $\mu +\sum t^{i}\varphi _{i}$ of $\mu $ is isomorphic to $\mu .$
By the commutativity of $\mu $ we can assume that every $\varphi _{i}$ is a
symmetric bilinear mapping.Then considering only the symmetric cochains $%
\varphi ,$ we can define the commutative Harrison cohomology of $\mu $ as follows : 
$$
Z_{s}^{2}(\mu ,\mu )=\{\varphi \in Sym(\mathbb{R}^{3}\times \mathbb{R}^{3},%
\mathbb{R}^{3})\quad \mid \quad \delta _{\mu }\varphi =0\}
$$
where 
$$
\delta _{\mu }\varphi (x,y,z)=\mu (x,\varphi (y,z))-\varphi (\mu (x,y),z)-\mu
(\varphi (x,y),z)+\varphi (x,\mu (y,z)).
$$
We can note that for every $f\in End(\mathbb{R}^{3})$ the coboundary $\delta
_{\mu }f\in Sym(\mathbb{R}^{3}\times \mathbb{R}^{3},\mathbb{R}^{3})$ and
then $\delta _{\mu }f\in Z_{s}^{2}(\mu ,\mu ).$

\begin{proposition}
If $H_{s}^{2}(\mu ,\mu )=\frac{Z_{s}^{2}(\mu ,\mu )}{B_{s}^{2}(\mu ,\mu )}%
=\{0\}$, the corresponding affine structure on $\mathbb{R}^{3}$ is 
rigid.
\end{proposition}

\section{Affine structures on the 2-dimensional abelian Lie algebra}

\subsection{ Classification of commutative associative $2$-dimensional real algebras}

$\bullet $ Suppose that the algebra $A$ is unitary. Then its
law is isomorphic to 
$$
\left\{ 
\begin{array}{l}
\mu _{1}(e_{1},e_{1})=e_{1} \\ 
\mu _{1}(e_{1},e_{2})=e_{2} \\ 
\mu _{1}(e_{2},e_{1})=e_{2} \\ 
\mu _{1}(e_{2},e_{2})=e_{2}
\end{array}
\right. \qquad ;\qquad \left\{ 
\begin{array}{l}
\mu _{2}(e_{1},e_{1})=e_{1} \\ 
\mu _{2}(e_{1},e_{2})=e_{2} \\
\mu _{2}(e_{2},e_{1})=e_{2} \\  
\mu _{2}(e_{2},e_{2})=0
\end{array}
\right. \qquad ;\qquad \left\{ 
\begin{array}{l}
\mu _{3}(e_{1},e_{1})=e_{1} \\ 
\mu _{3}(e_{1},e_{2})=e_{2} \\
\mu _{3}(e_{2},e_{1})=e_{2} \\  
\mu _{3}(e_{2},e_{2})=-e_{1}
\end{array}
\right.
$$

$\bullet $ Suppose that $A$ is nilpotent (and not unitary). The law is
isomorphic to 
$$
\mu _{4}(e_{1},e_{1})=e_{2} \qquad ;\qquad \mu _{5}=0
$$

$\bullet $ Suppose that $A$ is non nilpotent and non unitary. Then its law
is isomorphic to 
$$
\mu _{6}(e_{1},e_{1})=e_{1}.
$$

\subsection{Description of the affine structures}

\begin{proposition}
There are 6 affinely non-equivalent affine structures on the 2-dimensional abelian
Lie algebra.
\end{proposition}

In the following table we give the affine structures on the 2-dimensional
Lie algebra, the corresponding action and precise the completeness or not of
these structures. 
$$
\begin{tabular}{|l|l|l|l|l|}
\hline
&  & {\small affine structure} & {\small affine action} & {\small complete}
\\ \hline
$A_{1}$ &  & $\left( 
\begin{array}{ccc}
a & 0 & a \\ 
b & a+b & b \\ 
0 & 0 & 0
\end{array}
\right) $ & $\left( 
\begin{array}{ccc}
e^{a} & 0 & e^{a}-1 \\ 
e^{a}(e^{b}-1) & e^{a}e^{b} & e^{a}(e^{b}-1) \\ 
0 & 0 & 1
\end{array}
\right) $ & {\small NO} \\ \hline
$A_{2}$ & $\quad $ & $\left( 
\begin{array}{ccc}
a & 0 & a \\ 
b & a & b \\ 
0 & 0 & 0
\end{array}
\right) $ & $\left( 
\begin{array}{ccc}
e^{a} & 0 & e^{a}-1 \\ 
be^{a} & e^{a} & be^{a} \\ 
0 & 0 & 1
\end{array}
\right) $ & {\small NO} \\ \hline
$A_{3}$ &  & $\left( 
\begin{array}{ccc}
a & -b & a \\ 
b & a & b \\ 
0 & 0 & 0
\end{array}
\right) $ & $\left( 
\begin{array}{ccc}
e^{a}\cos b & -e^{a}\sin b & 1-e^{a}\cos b \\ 
e^{a}\sin b & e^{a}\cos b & e^{a}\sin b \\ 
0 & 0 & 1
\end{array}
\right) $ & {\small NO} \\ \hline
$A_{4}$ & $\quad $ & $\left( 
\begin{array}{ccc}
0 & 0 & a \\ 
a & 0 & b \\ 
0 & 0 & 0
\end{array}
\right) $ & $\left( 
\begin{array}{ccc}
1 & 0 & a \\ 
a & 1 & \frac{a^{2}}{2}+b \\ 
0 & 0 & 1
\end{array}
\right) $ & {\small YES} \\ \hline
$A_{5}$ & $\quad $ & $\left( 
\begin{array}{ccc}
0 & 0 & a \\ 
0 & 0 & b \\ 
0 & 0 & 0
\end{array}
\right) $ & $\left( 
\begin{array}{ccc}
1 & 0 & a \\ 
0 & 1 & b \\ 
0 & 0 & 1
\end{array}
\right) $ & {\small YES} \\ \hline
$A_{6}$ &  & $\left( 
\begin{array}{ccc}
a & 0 & a \\ 
0 & 0 & b \\ 
0 & 0 & 0
\end{array}
\right) $ & $\left( 
\begin{array}{ccc}
e^{a} & 0 & e^{a}-1 \\ 
0 & 1 & b \\ 
0 & 0 & 1
\end{array}
\right) $ & {\small NO} \\ \hline
\end{tabular}
$$

\subsection{\protect\bigskip On the group of affine transformations}

To each affine structure $A_{i}$ corresponds a flat affine connection
without torsion $\nabla ^{i}$ on the abelian Lie group $G$.

The set of all affine transformations of $(G,\nabla ^{i})$ is a Lie group,
noted $Aff(G,\nabla ^{i}).$ Its Lie algebra is the set of complete affine
vector fields ([3]), that is, complete vector field satisfying 
$$
\left[ X,\nabla ^{i}_{Y}Z\right] =\nabla ^{i}_{\left[ X,Y\right] }Z+\nabla ^{i}_{Y}%
\left[ X,Z\right]
$$
If $\nabla ^{i}$ is complete (in our case $i=4,5$ ) then the Lie algebra of $%
Aff(G,\nabla ^{i})$ is the Lie algebra of affine vector fields $aff(G,\nabla
^{i})$ and $Aff(G,\nabla ^{i})$ acts transitively on $G$.

Let us consider the corresponding affine action to $\nabla _{i}$ described
in the previous section. The translation part defines an open set $U_{i}\subset 
\mathbb{R}^{2}$. For the complete case we obtain $U_{i}=\mathbb{R}^{2}$. For
the non
 complete case we have 
\begin{eqnarray*}
U_{1} &=&U_{2}=U_{6}=\left\{ \left( x,y\right) \in \mathbb{R}%
^{2}\,/\,x>-1\right\} \\
U_{3} &=&\left\{ \left( x,y\right) \in \mathbb{R}^{2}\,,\,\left( x,y\right)
\neq \left( 1,0\right) \right\}
\end{eqnarray*}
If $\phi $ is an affine transformation which leaves $U_{i}$ invariant, for $%
i=1,2,6$ the matrix of $\phi $ is seen to have the following form 
$$
\left( 
\begin{array}{lll}
a & 0 & e^{t}-1 \\ 
b & c & u \\ 
0 & 0 & 1
\end{array}
\right) .
$$
The group $Aff(G,\nabla ^{i})$ is the maximal group included in the semi
group generated by the previous matrices.

Then the group $Aff(G,\nabla ^{i})=B_{2}\times \mathbb{R}^2$, where $B_{2}$ is
the subgroup of $GL(2,\mathbb{R)}$ constituted of triangular matrices.

\subsection{Rigid affine structures on the 2-dimensional abelian Lie algebra}

\begin{lemma}
Every infinitesimal deformation of an unitary associative algebra in $%
\mathcal{A}^{c}(n)$ is unitary.
\end{lemma}

The proof is based on the study of perturbations of idempotent elements made
in [11]. Let $X\,$be an idempotent element in an associative algebra of law $\mu .$
The operators 
$$
l_{X} : Y{\longrightarrow }XY
$$
and 
$$
r_{X} : Y {\longrightarrow }YX
$$
are simultaneously diagonalizable and the eigenvalues are respectively
(1,...,1,0,  ...,0) and (1,1,1,...,1,0,...,0). This set of eigenvalues is
called bisystem associated to $X.$ If $X$ corresponds to the identity, the
bisystem is $\left\{ \left( 1,...,1\right) \left( 1,...,1\right) \right\} .$

From [11], if $\mu ^{\prime }$ is a perturbation of $\mu ,$ there exists in $%
\mu ^{\prime }$ an idempotent element $X^{\prime }$ such that $b(X^{\prime
})=b(X).$ Consider the perturbation of the identity in $\mu ^{\prime }.$ As
the bisystem corresponds to $\left\{ \left\{ 1,...,1\right\} ,\left\{
1,...,1\right\} \right\} $, we can conclude that $X^{\prime }$ is the
identity of $\mu ^{\prime }.$

\medskip

\noindent \textbf{Consequence. }The set of unitary associative algebra is open
in $\mathcal{A}^{c}(n)$.

\medskip

Let us note by $\theta (\mu )$ the orbit of the law $\mu $ in $\mathcal{A}%
^{c}(n).$ From the previous classification, we see that 
$$
\mu _{i}\in \overline{\theta (\mu _{1})}
$$
for $i=2,5,6$. Then we have

\begin{theorem}
The affine structures $A_{1}$ and $A_{3}$ on the 2-dimensional abelian Lie
algebra are rigid. The other structures can be deformed into $A_{1}$ or $%
A_{3}.$
\end{theorem}

\noindent \textbf{Consequence. }Not any complete affine structure is rigid.

\subsection{Invariant affine structure on $\mathbb{T}^{2}$}

The compact abelian Lie group $\mathbb{T}^{2}$ is defined by $\mathbb{T}^{2}=%
\mathbb{R}^{2}/\mathbb{Z}^{2}$ identifying $(x,y)$ with $(x+p,y+q)$, $%
(p,q)\in \mathbb{Z}^{2}.$

\begin{proposition}
Only the structures $A_{4}$ and $A_{5}$ induce affine structures on the torus 
$\mathbb{T}^{2}$.
\end{proposition}

\noindent {\it Proof.} It is easy to see that the affine action associated to $%
A_{1},A_{2},A_{3}$ and $A_{6}$ are incompatible with the lattice defined by $%
\mathbb{Z}^{2}$. Thus only complete structures provide affine structures on 
$\mathbb{T}^{2}$. For $A_{4}$ we obtain the following affine transformations 
$$
\left( x,y\right) \rightarrow \left( x+p,px+y+\left( q+p^{2}\right) \right) 
$$
This structure on $\mathbb{T}^{2}$ is not euclidean. For $A_{5}$ the affine
structure on $\mathbb{T}^{2}$ which is euclidean corresponds to the
transformations 
$$
\left( x,y\right) \rightarrow \left( x+p,y+q\right) 
$$

\begin{remark}
In this proposition we find again, for the particular case of the torus, a
classical result of Kuiper \mbox{[9]} giving the classification of affine
structures on surfaces.
\end{remark}

\section{Affine structures on the 3-dimensional abelian Lie algebra}

\subsection{Classification of 3-dimensional commutative associative real algebras}

Let us begin by describing the classification of real associative
commutative algebras. Let $\frak{a}$ be a 3-dimensional (not necessarily
unitary) real commutative associative algebra.

If $\frak{a}$ is simple, then $\frak{a}$ is, following Wedderburn's theorem,
isomorphic to $\left( M_{1}(\mathbb{R)}\right) ^{3}$, $M_{1}(\mathbb{%
R)\oplus }M_{1}(\mathbb{C)}$, where $M_{n}(D\mathbb{)}$ is a matrix algebra
on a division algebra on $\mathbb{R}$, that is $D=\mathbb{R}$ or $\mathbb{C}$%
. This gives the following algebras 
$$
\left\{ 
\begin{array}{l}
\mu _{1}(e_{1},e_{i})=e_{i}\qquad i=1,2,3 \\ 
\mu _{1}(e_{i},e_{1})=e_{i}\qquad i=1,2,3 \\ 
\mu _{1}(e_{2},e_{2})=e_{2} \\ 
\mu _{1}(e_{3},e_{3})=e_{3}
\end{array}
\right. \qquad ,\qquad \left\{ 
\begin{array}{l}
\mu _{2}(e_{1},e_{i})=e_{i}\qquad i=1,2,3 \\ 
\mu _{2}(e_{i},e_{1})=e_{i}\qquad i=1,2,3 \\ 
\mu _{2}(e_{2},e_{2})=e_{2} \\ 
\mu _{2}(e_{3},e_{3})=e_{2}-e_{1}
\end{array}
\right.
$$

If $\frak{a}$ is not simple, then $\frak{a}=J(\frak{a})\oplus \frak{s}$,
where $\frak{s}$ is simple and $J(\frak{a})$ is the Jacobson radical of $\frak{%
a.}$ If \thinspace $\frak{s}=\left( M_{1}(\mathbb{R)}\right) ^{2}$, we
obtain 
$$
\left\{ 
\begin{array}{l}
\mu _{3}(e_{1},e_{i})=e_{i}\qquad i=1,2,3 \\ 
\mu _{3}(e_{i},e_{1})=e_{i}\qquad i=1,2,3 \\ 
\mu _{3}(e_{2},e_{2})=e_{2}
\end{array}
\right.
$$
If \thinspace $\frak{s}=M_{1}(\mathbb{C)}$%
$$
\left\{ 
\begin{array}{l}
\mu _{4}(e_{1},e_{i})=e_{i}\qquad i=1,2,3 \\ 
\mu _{4}(e_{i},e_{1})=e_{i}\qquad i=1,2,3 \\ 
\mu _{4}(e_{3},e_{3})=e_{2}
\end{array}
\right.
$$
where $J(\frak{a})$ is not abelian or 
$$
\left\{
\begin{array}{l}
\mu _{5}(e_{1},e_{i})=e_{i}\qquad i=1,2,3  \\
\mu _{5}(e_{i},e_{1})=e_{i}\qquad i=1,2,3
\end{array}
\right.
$$
where $J(\frak{a})$ is abelian.

Suppose that $\frak{a}$ is not unitary. As the Levi decomposition also
holds in this case, we have the following possibilities : $\frak{a}\simeq
\left( M_{1}(\mathbb{R)}\right) ^{2}\oplus J(\frak{a})$ or $M_{1}(\mathbb{C)}%
\oplus J(\frak{a})$. This gives the following algebras 
$$
\left\{ 
\begin{array}{c}
\mu _{6}(e_{1},e_{1})=e_{1} \\ 
\mu _{6}(e_{2},e_{2})=e_{2}
\end{array}
\right. \qquad ,\qquad \left\{ 
\begin{array}{l}
\mu _{7}(e_{1},e_{1})=e_{1} \\ 
\mu _{7}(e_{1},e_{2})=e_{2} \\ 
\mu _{7}(e_{2},e_{1})=e_{2} \\ 
\mu _{7}(e_{2},e_{2})=-e_{1}
\end{array}
\right.
$$
$$
\left\{ \mu _{8}(e_{1},e_{1})=e_{1}\qquad \right. \qquad ,\qquad \left\{ 
\begin{array}{c}
\mu _{9}(e_{1},e_{1})=e_{1} \\ 
\mu _{9}(e_{1},e_{2})=e_{2} \\
\mu _{9}(e_{2},e_{1})=e_{2}
\end{array}
\right.
$$
$$
\left\{ 
\begin{array}{c}
\mu _{10}(e_{1},e_{1})=e_{1} \\ 
\mu _{10}(e_{2},e_{2})=e_{3}
\end{array}
\right. \qquad
$$
If moreover $\frak{a}$ is nilpotent, then it is isomorphic to one of the following a$\lg $%
ebras 
$$
\left\{ 
\begin{array}{c}
\mu _{11}(e_{1},e_{1})=e_{2} \\ 
\mu _{11}(e_{3},e_{3})=e_{2}
\end{array}
\right. \qquad ,\qquad \left\{ 
\begin{array}{l}
\mu _{12}(e_{1},e_{1})=e_{2} \\ 
\mu _{12}(e_{3},e_{3})=-e_{2}
\end{array}
\right.
$$
$$
\left\{ \mu _{13}(e_{1},e_{1})=e_{2}\right. \qquad ,\qquad \left\{ 
\begin{array}{l}
\mu _{14}(e_{1},e_{1})=e_{2} \\ 
\mu _{14}(e_{1},e_{2})=e_{3} \\
\mu _{14}(e_{2},e_{1})=e_{3}
\end{array}
\right.
$$
$$
\left\{ \mu _{15}(e_{i},e_{j})=0\qquad i,j\in \{1,2,3\}\right.
$$

\begin{theorem}
Every 3-dimensional real commutative associative Lie algebra $\frak{a}$ is
isomorphic to one of the algebras $\frak{a}_{i},$ $i=1,2,...,15.$

\noindent If $\frak{a}$ is nilpotent, $\frak{a}$ $\ $is isomorphic to $\frak{a}_{i},$ $%
i=11,...,15.$
\end{theorem}

\subsection{Affine structures on $\mathbb{R}^{3}$}

For each associative algebra of the previous classification we can describe
the affine action on $\mathbb{R}^{3}$.

\begin{theorem}
There exist 15 invariant affinely non-equivalent affine structures on the
3-dimensional abelian Lie algebra.\ They are given by : 
\end{theorem}

\bigskip

$$
\begin{tabular}{|l|l|}
\hline
& \\
$A_{1}$ & $(x,y,z)\rightarrow \left\{ 
\begin{array}{l}
e^{a}x+e^{a}-1, \\ 
e^{a}\left( e^{b}-1\right) x+e^{a+b}y +e^a(e^b-1), \\ 
e^{a}\left( e^{c}-1\right) x+e^{a+c}z+e^{a}\left( e^{c}-1\right)
\end{array}
\right. $ \\ \hline
& \\
$A_{2}$ & $(x,y,z)\rightarrow \left\{ 
\begin{array}{l}
xe^{a}\cos c-ze^{a}\sin c+e^{a}\cos c-1, \\ 
\left( -e^{a}\cos c+e^{a+b}\right) x+e^{a+b}y+ze^{a}\sin c-e^{a}\cos
c+e^{a+b}, \\ 
\ xe^{a}\sin c+ze^{a}\cos c+e^{a}\sin c
\end{array}
\right. $ \\ \hline
& \\
$A_{3}$ & $(x,y,z)\rightarrow \left\{ 
\begin{array}{l}
e^{a}x+e^{a}-1, \\ 
e^{a}\left( e^{b}-1\right) x+e^{a+b}y+e^{a}\left( e^{b}-1\right) , \\ 
ce^{a} x+e^{a}z+ce^{a}
\end{array}
\right. $ \\ \hline
& \\
$A_{4}$ & $(x,y,z)\rightarrow \left\{ 
\begin{array}{l}
e^{a}x-1+e^{a}, \\ 
(b+\frac{c^2}{2})e^ax+e^{a}y+ce^az+(b+\frac{c^2}{2})e^a, \\ 
\ ce^ax+e^{a}z+ce^a
\end{array}
\right. $ \\ \hline
& \\
$A_{5}$ & $(x,y,z)\rightarrow \left\{ 
\begin{array}{l}
e^{a}x-1+e^{a}, \\ 
be^{a}x+e^{a}y+be^{a}, \\ 
ce^{a}x+e^{a}z+ce^{a}
\end{array}
\right. $ \\ \hline
& \\
$A_{6}$ & $(x,y,z)\rightarrow \left\{ 
\begin{array}{l}
e^{a}x-1+e^{a}, \\ 
e^{b}y-1+e^{b}, \\ 
z+c
\end{array}
\right. $ \\ \hline
& \\
$A_{7}$ & $(x,y,z)\rightarrow \left\{ 
\begin{array}{l}
xe^{a}\cos b-ye^{a}\sin b-1+e^{a}\cos b, \\ 
xe^{a}\sin b+ye^{a}\cos b+e^{a}\sin b, \\ 
z+c
\end{array}
\right. $ \\ \hline
& \\
$A_{8}$ & $(x,y,z)\rightarrow \left\{ 
\begin{array}{l}
e^{a}x-1+e^{a}, \\ 
y+b, \\ 
z+c
\end{array}
\right. $ \\ \hline
& \\
$A_{9}$ & $(x,y,z)\rightarrow \left\{ 
\begin{array}{l}
e^{a}x-1+e^{a}, \\ 
be^{a}x+e^{a}y+be^{a}, \\ 
z+c
\end{array}
\right. $ \\ \hline
& \\
$A_{10}$ & $(x,y,z)\rightarrow \left\{ 
\begin{array}{l}
e^{a}x-1+e^{a}, \\ 
y+b, \\ 
by+z+\frac{b^2}{2}c
\end{array}
\right. $ \\ \hline
& \\
$A_{11}$ & $(x,y,z)\rightarrow \left\{ 
\begin{array}{l}
x+a, \\ 
ax+y+cz+b+\frac{1}{2}(a^2+c^2), \\ 
z+c
\end{array}
\right. $ \\ \hline
& \\
$A_{12}$ & $(x,y,z)\rightarrow \left\{ 
\begin{array}{l}
x+a, \\ 
ax+y-cz+b+\frac{1}{2}(a^2-c^2), \\ 
z+c
\end{array}
\right. $ \\ \hline
\end{tabular}
$$
$$
\begin{tabular}{|l|l|}
\hline
& \\
$A_{13}$ & $(x,y,z)\rightarrow \left\{ 
\begin{array}{l}
x+a, \\ 
ax+y+b+\frac{a^2}{2}, \\ 
z+c
\end{array}
\right. $ \\ 
& \\
\hline
& \\
$A_{14}$ & $(x,y,z)\rightarrow \left\{ 
\begin{array}{l}
x+a, \\ 
ax+y+b+\frac{a^2}{2}, \\ 
(b+\frac{a^2}{2})x+ay+z+\frac{a^3}{6}+ab+c
\end{array}
\right. $ \\ 
& \\
\hline
& \\
$A_{15}$ & $(x,y,z)\rightarrow \left\{ 
\begin{array}{l}
x+a, \\ 
y+b, \\ 
z+c
\end{array}
\right. $ \\ 
& \\
\hline
\end{tabular}
$$
{\it where $a,b,c \in \Bbb{R}.$ Only the structures $A_{i},$ $i=11,..,15$ are complete.}

\subsection{Rigid affine structures}

Recall that an affine structure on $\mathbb{R}^{3}$ is called rigid if the
corresponding commutative and associative real algebra is rigid in the
algebraic variety $\mathcal{A}^{c}(3).$

As the laws $\mu _{1}$ and $\mu _{2}$ are semi-simple associative,
their second cohomological group is trivial.\ These structures are rigid.

Consider the remaining laws $\mu _{i}.$ We can easily compute the linear space $%
H_{s}^{2}(\mu ,\mu )$ and presente the results in the following tables :

\bigskip Unitary case :

$$
\begin{tabular}{|l|l|l|}
\hline
& & \\
laws & $dimH_{s}^{2}(\mu ,\mu )$ & basis of $H_{s}^{2}(\mu ,\mu )$ \\ 
& & \\
\hline
& & \\
$\mu _{3}$ & $1$ & $\varphi (e_{3},e_{3})=e_{2}-e_{1}$ \\ 
& & \\
\hline
& & \\
$\mu _{4}$ & $2$ & $\left\{ 
\begin{array}{c}
\varphi _{1}(e_{2},e_{2})=e_{2} \\ 
\varphi _{2}(e_{2},e_{3})=e_{3.}
\end{array}
\right. $ \\ 
& & \\
\hline
& & \\
$\mu _{5}$ & $4$ & $\left\{ 
\begin{array}{c}
\varphi _{1}(e_{2},e_{2})=e_{2} \\ 
\varphi _{2}(e_{2},e_{3})=e_{3.}
\end{array}
\right. ;\left\{ 
\begin{array}{c}
\varphi _{3}(e_{3},e_{3})=e_{2} \\ 
\varphi _{4}(e_{3},e_{3})=e_{3.}
\end{array}
\right. $ \\ 
& & \\
\hline
\end{tabular}
$$

\bigskip Non unitary case : 
$$
\begin{tabular}{|l|l|l|}
\hline
& & \\
laws & $dimH_{s}^{2}$ & basis of $H_{s}^{2}(\mu ,\mu )$ \\ 
& & \\ \hline
& & \\
$\mu _{6}$ & $1$ & $\varphi _{1}(e_{3},e_{3})=e_{3}$ \\
& & \\
\hline
& & \\
$\mu _{7}$ & $1$ & $\varphi _{1}(e_{3},e_{3})=e_{3}$ \\ 
& & \\
\hline
& & \\
$\mu _{8}$ & $6$ & $\left\{ 
\begin{array}{c}
\varphi _{1}(e_{2},e_{2})=e_{2} \\ 
\varphi _{2}(e_{2},e_{2})=e_{3.}
\end{array}
\right. ;\left\{ 
\begin{array}{c}
\varphi _{3}(e_{3},e_{3})=e_{2} \\ 
\varphi _{4}(e_{3},e_{3})=e_{3.}
\end{array}
\right. ;\left\{ 
\begin{array}{c}
\varphi _{5}(e_{2},e_{3})=e_{2} \\ 
\varphi _{6}(e_{3},e_{3})=e_{3.}
\end{array}
\right. $ \\
& & \\
 \hline
& & \\
$\mu _{9}$ & $3$ & $\left\{ 
\begin{array}{c}
\varphi _{1}(e_{2},e_{2})=e_{1} \\ 
\varphi _{2}(e_{2},e_{3})=e_{2.}
\end{array}
\right. ;\quad \varphi _{3}(e_{3},e_{3})=e_{3}$ \\ 
& & \\
\hline
& & \\
$\mu _{10}$ & $1$ & $\varphi _{1}(e_{2},e_{3})=e_{2},\varphi
_{1}(e_{3},e_{3})=e_{3}.$ \\ 
& & \\
\hline
\end{tabular}
$$

Nilpotent and complete case : 
$$
\begin{tabular}{|l|l|l|}
\hline
& & \\
laws & $dimH_{s}^{2}$ & basis of $H_{s}^{2}(\mu ,\mu )$ \\ 
& & \\
\hline
& & \\
$\mu _{11}$ & $3$ & $\left\{ 
\begin{array}{c}
\varphi _{1}(e_{1},e_{3})=e_{1} \\ 
\varphi _{1}(e_{2},e_{3})=e_{2.}
\end{array}
\right. ;\left\{ 
\begin{array}{c}
\varphi _{2}(e_{1},e_{3})=e_{3} \\ 
\varphi _{2}(e_{3},e_{3})=e_{1.}
\end{array}
\right. ;\left\{ 
\begin{array}{c}
\varphi _{3}(e_{3},e_{3})=e_{3}.
\end{array}
\right. $ \\ 
& & \\
\hline
$\mu _{12}$ & $4$ & $
\begin{array}{c}
\varphi _{1}(e_{1},e_{2})=e_{2}. \\ 
\varphi _{2}(e_{3},e_{3})=e_{3.}
\end{array}
;\left\{ 
\begin{array}{c}
\varphi _{3}(e_{1},e_{3})=e_{3} \\ 
\varphi _{3}(e_{3},e_{3})=e_{1.}
\end{array}
\right. ;\left\{ 
\begin{array}{l}
\varphi _{4}(e_{1},e_{1})=e_{3} \\ 
\varphi _{4}(e_{1},e_{3})=e_{1} \\ 
\varphi _{4}(e_{2},e_{3})=2e_{2}.
\end{array}
\right. $ \\ 
& & \\
\hline
$\mu _{13}$ & $7$ & $
\begin{array}{c}
\left\{ 
\begin{array}{c}
\varphi _{1}(e_{1},e_{2})=e_{1} \\ 
\varphi _{1}(e_{2},e_{2})=e_{2.}
\end{array}
\right. \\ 
\varphi _{2}(e_{1},e_{2})=e_{2}.
\end{array}
; 
\begin{array}{c}
\varphi _{3}(e_{1},e_{2})=e_{3} \\ 
\left\{ 
\begin{array}{c}
\varphi _{4}(e_{1},e_{3})=e_{1} \\ 
\varphi _{4}(e_{2},e_{3})=e_{2}.
\end{array}
\right.
\end{array}
\begin{array}{c}
\varphi _{5}(e_{1},e_{3})=e_{3}. \\ 
\varphi _{6}(e_{3},e_{3})=e_{2}. \\ 
\varphi _{7}(e_{3},e_{3})=e_{3}.
\end{array}
$ \\ 
\hline
& & \\
$\mu _{14}$ & $3$ & $\left\{ 
\begin{array}{c}
\varphi _{1}(e_{1},e_{3})=e_{1} \\ 
\varphi _{1}(e_{2},e_{2})=e_{1} \\ 
\varphi _{1}(e_{2},e_{3})=e_{2} \\ 
\varphi _{1}(e_{3},e_{3})=e_{3}.
\end{array}
\right. ;\left\{ 
\begin{array}{c}
\varphi _{2}(e_{1},e_{3})=e_{2} \\ 
\varphi _{2}(e_{2},e_{2})=e_{2} \\ 
\varphi _{2}(e_{2},e_{3})=e_{3}.
\end{array}
\right. ;\left\{ 
\begin{array}{c}
\varphi _{3}(e_{1},e_{3})=e_{3} \\ 
\varphi _{3}(e_{2},e_{2})=e_{3.}
\end{array}
\right. $ \\ 
\hline
& & \\
$\mu _{15}$ & $18$ &  \\ 
\hline
\end{tabular}
$$
In this box we suppose that $\varphi (e_{i},e_{j})=\varphi (e_{j},e_{i})$
and the non defined $\varphi (e_{s},e_{t})$ are equal to zero.

We will note by $\mu _{i}\rightarrow \mu _{j}$ when $\mu _{i}\in \overline{%
\mathcal{O(}\mu _{j})}.$ We obtain the following diagramm 
\[
\begin{tabular}{lllllllllllll}
&  &  &  &  &  &  &  &  &  &  &  &  \\ 
&  &  & $\mu _{12}$ &  &  &  & $\mu _{14}$ &  &  &  &  &  \\ 
&  &  &  & $\searrow $ &  & $\swarrow $ &  &  &  &  &  &  \\ 
& $\mu _{11}$ & $\rightarrow $ & $\mu _{10}$ & $\rightarrow $ & $\mu _{1}$ & 
$\leftarrow $ & $\mu _{6}$ & $\leftarrow $ & $\mu _{8}$ & $\longrightarrow $
& $\mu _{9}$ &  \\ 
& $\uparrow $ &  &  &  & $\uparrow $ &  &  &  &  &  & $\downarrow $ &  \\ 
& $\mu _{13}$ &  &  &  & $\mu _{4}$ & $\leftarrow $ & $\mu _{5}$ & $%
\rightarrow $ & $\mu _{3}$ & $\longleftarrow $ & $\mu _{7}$ &  \\ 
&  &  &  &  &  &  &  &  & $\downarrow $ &  &  &  \\ 
&  &  &  &  &  &  &  &  & $\mu _{2}$ &  &  &  \\ 
&  &  &  &  &  &  &  &  &  &  &  & 
\end{tabular}
\]

\begin{theorem}
There exist two rigid affine structures on the 3-dimensional abelian Lie
algebra. They are the structures associated to the semi simple associative
algebras $\mu _{1}$ and $\mu _{2}.$
\end{theorem}

\subsection{Invariant affine structures on $\mathbb{T}^{3}$}

We can see that the affine actions $A_{1}$ to $A_{10}$ are compatible with
the action of $\mathbb{Z}^{3}$ on $\mathbb{R}^{3}$ if the $\exp $onentials
which appear in the analytic expressions of the affine transformations are
equal to $1$. This gives the identity for $A_{1}$ and $A_{3}$. As $A_{2}$ is
incompatible with the action of $\mathbb{Z}^{3}$ for any values of the
parameters $a,b,c,$ the affine structures corresponding to the unitary cases
are given by $A_{4}$ and $A_{5}$ for $a=0$. This corresponds to the
following affine structures on the torus $\mathbb{T}^{3}:$%
$$
(x,y,z)\in \mathbb{R}^{3}/\mathbb{Z}^{3}\rightarrow \left\{ 
\begin{array}{l}
x \\ 
px+y+qz+p \\ 
qx+z+q
\end{array}
\right.
$$
and 
$$
(x,y,z)\in \mathbb{R}^{3}/\mathbb{Z}^{3}\rightarrow \left\{ 
\begin{array}{l}
x \\ 
px+y+p \\ 
qx+z+q
\end{array}
\right.
$$
with $p,q\in \mathbb{Z}$.

For the actions $A_{6}$ to $A_{10},$ they induce affine actions on the
torus if $a=b=0$ for $A_{6}$ and $a=0$ for $A_{7}$ to $A_{10}$ but this
appears as a particular case of $A_{11},A_{13}$ and $A_{14}.$ Let us
examine the complete and nilpotent cases; we find the following affine
structure on $\mathbb{T}^{3}:$%
$$
\begin{tabular}{|l|}
\hline
$(x,y,z)\in \mathbb{R}^{3}/\mathbb{Z}^{3}\rightarrow \left\{ 
\begin{array}{l}
x+p, \\ 
px+y+rz+q, \\ 
z+r
\end{array}
\right. $ \\ \hline
$(x,y,z)\in \mathbb{R}^{3}/\mathbb{Z}^{3}\rightarrow \left\{ 
\begin{array}{l}
x+p, \\ 
px+y-rz+q, \\ 
z+r
\end{array}
\right. $ \\ \hline
$(x,y,z)\in \mathbb{R}^{3}/\mathbb{Z}^{3}\rightarrow \left\{ 
\begin{array}{l}
x+p, \\ 
px+y+q, \\ 
z+r
\end{array}
\right. $ \\ \hline
$(x,y,z)\in \mathbb{R}^{3}/\mathbb{Z}^{3}\rightarrow \left\{ 
\begin{array}{l}
x+p, \\ 
px+y+q, \\ 
qx+py+z+r
\end{array}
\right. $ \\ \hline
$(x,y,z)\in \mathbb{R}^{3}/\mathbb{Z}^{3}\rightarrow \left\{ 
\begin{array}{l}
x+p, \\ 
y+q, \\ 
z+r
\end{array}
\right. $ \\ \hline
\end{tabular}
$$

\begin{theorem}
There exist 7 affine structures on the torus $\mathbb{T}^{3}.$ They
correspond to the following affine crystallographic subgroups of $Aff(%
\mathbb{R}^{3}):$%
$$
\begin{tabular}{|l|l|}
\hline
& \\
$\Gamma _{1}=\left\{ \left( 
\begin{array}{llll}
1 & 0 & 0 & 0 \\ 
p & 1 & q & p \\ 
q & 0 & 1 & q \\ 
0 & 0 & 0 & 1
\end{array}
\right) \right\} $ & $\Gamma _{2}=\left\{ \left( 
\begin{array}{llll}
1 & 0 & 0 & 0 \\ 
p & 1 & 0 & p \\ 
q & 0 & 1 & q \\ 
0 & 0 & 0 & 1
\end{array}
\right) \right\} $ \\ 
& \\
\hline
& \\
$\Gamma _{3}=\left\{ \left( 
\begin{array}{llll}
1 & 0 & 0 & p \\ 
p & 1 & r & q \\ 
0 & 0 & 1 & r \\ 
0 & 0 & 0 & 1
\end{array}
\right) \right\} $ & $\Gamma _{4}=\left\{ \left( 
\begin{array}{llll}
1 & 0 & 0 & p \\ 
p & 1 & -r & q \\ 
0 & 0 & 1 & r \\ 
0 & 0 & 0 & 1
\end{array}
\right) \right\} $ \\ 
& \\
\hline
& \\
$\Gamma _{5}=\left\{ \left( 
\begin{array}{llll}
1 & 0 & 0 & p \\ 
p & 1 & 0 & q \\ 
0 & 0 & 1 & r \\ 
0 & 0 & 0 & 1
\end{array}
\right) \right\} $ & $\Gamma _{6}=\left\{ \left( 
\begin{array}{llll}
1 & 0 & 0 & p \\ 
p & 1 & 0 & q \\ 
q & p & 1 & r \\ 
0 & 0 & 0 & 1
\end{array}
\right) \right\} $ \\ 
& \\
\hline
& \\
$\Gamma _{7}=\left\{ \left( 
\begin{array}{llll}
1 & 0 & 0 & p \\ 
0 & 1 & 0 & q \\ 
0 & 0 & 1 & r \\ 
0 & 0 & 0 & 1
\end{array}
\right) \right\} $ &  \\ 
& \\
\hline
\end{tabular}
$$
\end{theorem}

\begin{remark}
The (complete) nilpotent and unimodular cases are completely classified in
[6].\ 
\end{remark}

\bigskip

\noindent REFERENCE

\noindent [1] Ancochea Bermudez J.M., \textit{On the rigidity of solvable Lie
algebras}. ASI NATO, Serie C247, 403-445, 1986

\noindent [2] Ancochea Bermudez J.M., Campoamor Stursberg O.R.,
\textit{On Lie algebras whose nilradical is (n-p)-filiform}. Comm. in Algebra,
\textbf{29 (1)}, 427-450, 2001  

\noindent [3] Bordemann M, Medina A., \textit{Le groupe des transformations affines d'un groupe
de Lie \`a structure affine bi-invariante}. S\'eminaire G. Darboux. Univ. Montpellier II. 183-210. 1995-1996

\noindent [4] Burde D., \textit{Affine structures
on nilmanifolds}. Inter.
 J. Math., \textbf{7}, 599-616, 1996

\noindent [5] Dekimpe K., Ongenae V., \textit{On the number of abelian
left symmetric algebras. }P.A.M.S., \textbf{128 (11),} 3191-3200, 2000

\noindent [6] Fried D., Goldman W., \textit{Three dimensional affine
crystallographic groups.} Adv. Math\textit{.}\textbf{\ 47}, 1-49. 1983

\noindent [7] Gabriel P., \textit{Finite representation type is open.} L.N.
488, 1975

\noindent [8] Kim H., \textit{Complete left-invariant affine
structures on nilpotent Lie groups}. J. Diff. Geom\textit{.}\textbf{24, }%
373-394, 1986

\noindent [9] Kuiper N., \textit{Sur les surfaces localement affines}.
Colloque G\'{e}%
 om\'{e}trie diff\'{e}rentielle Strasbourg 79-87, 1953

\noindent [10] Mazzola G., \textit{The algebraic and geometric classification
of associative algebras of dimension five}. Manuscripta math., \textbf{27,}
81-101, 1979.

\noindent [11] Makhlouf N., \textit{The irreducible components of nilpotent
associative algebras}. Revista Matematica Univ. Complutense, \textbf{6 (1)},
27-40, 1993

\noindent [12] Nagano T., Yagi K., \textit{The affine structures on the
real two torus.}\ Osaka J.\ Math\textit{. }\textbf{11}, 181-210. 1974

\noindent [13] Remm E.,Goze M.,  \textit{Non complete affine structures on
Lie algebras of maximal class.} Intern. Journal of Math. and Math. Sc. \textbf{29 (2)}, 71-78, 2002

\noindent [14] Remm E. \textit{Structures affines sur les alg\`ebres de Lie et op\'erades Lie-admissibles.}
Th\`ese Universit\'e de Mulhouse. 2001

\bigskip 

\noindent ADRESS : 

\noindent Laboratoire de Math\'ematiques et Applications. Facult\'e des Sciences et Techniques, 4 rue des Fr\`eres Lumi\`ere,
F. 68093 MULHOUSE.

\noindent EMAIL : E.Remm@uha.fr ; M.Goze@uha.fr 
\end{document}